\documentclass[a4paper,12pt]{article}
\usepackage{amsmath, amssymb, amsthm, latexsym, amsopn, amsfonts}
\usepackage{amscd, stmaryrd}
\usepackage[cp1251]{inputenc}
\usepackage{mathtext}
\usepackage[english]{babel}
\usepackage{indentfirst}
\usepackage{color}
\def\{{\protect\lbrace}
\def\}{\protect\rbrace}
\def\S{$\mathchar"278$}

\begin{document}
\thispagestyle{empty}

\mbox{}

\begin{center}

{\bf A. A. Tuganbaev

$\aleph_0$-DISTRIBUTIVE MODULES AND RINGS}

\hfill \textbf{A.A. Tuganbaev}

\hfill National Research University "MPEI", Moscow, Russia

\hfill Lomonosov Moscow State University, Moscow, Russia
 
\hfill tuganbaev@gmail.com

\end{center}

\textbf{Abstract.}  Let $A$ be a ring with minimum condition on principal right ideals. It is proved that $\aleph_0$-distributive right (left) $A$-modules coincide with Artinian (Noetherian) right (left) $A$-modules. Rings, over which all right modules are direct sums of $\aleph_0$-distributive modules coincide with rings of finite representation type. Rings, whose right modules are semidistributive, coincide with Kawada rings, over basis rings of which all right modules are completely cyclic. The studies of Tuganbaev are supported by Russian Scientific Foundation, project 22-11-00052.

\textbf{Key words.} $\aleph_0$-distributive module, semidistributive module, distributive module, Koethe ring, Kawada ring, ring of finite representation type.

\textbf{MSC Classification. 16G60; 16D90; 16G60}

\section{Introduction}\label{section1}

All rings are assumed to be associative and with non-zero identity element. Modules are assumed to be unitary and, unless otherwise specified, all modules are right modules. The words of type "$A$ is a right Artinian ring" mean that "the module $A_A$ is Artinian". The words of type "$A$ is an Artinian ring" mean that "the modules $A_A$ and $_AA$ are Artinian". 

A module is said to be \textsf{$\aleph_0$-distributive} if it does not have subfactors which are infinite direct sums of isomorphic simple modules. We note that any distributive module (i.e., a module with distributive submodule lattice) is $\aleph_0$-distributive, since every distributive module does not have subfactors which are a direct sum of two isomorphic non-zero modules, \cite{Ste74}. The ring $\mathbb{Z}[x]$ is an example of a commutative $\aleph_0$-distributive ring which is not a distributive ring. 

\textbf{Remark 1.1.} It is clear that if all factor modules of the module $M$ do not contain infinite direct sums of non-zero submodules  (i.e. are \textsf{finite-dimensional}), such modules are called \textsf{quotient finite-dimensional} modules or \textsf{q.f.d.} modules), then the module $M$ is $\aleph_0$-distributive. Therefore, every Noetherian, or Artinian, or uniserial module is quotient finite-dimensional and, in particular, an $\aleph_0$-distributive module. (A module is said to be \textsf{uniserial}\label{uniserMod} if any two its (equivalently, cyclic) submodules are comparable with respect to inclusion. A direct sum of uniserial modules is called a \textsf{serial}\label{uniserMod} module.)\\
A direct product of an infinite number of fields is a commutative ($\aleph_0$-)distributive ring which is not quotient finite-dimensional. Every uniserial module is distributive. The ring of integers $\mathbb{Z}$ is a distributive non-uniserial $\mathbb{Z}$-module. 

\textbf{Remark 1.2.} A ring $A$ is said to be \textsf{right pure-semisimple} if $A$ is a right Artinian ring such that every right $A$-module is a direct sum of finitely generated modules.

A ring $A$ is called a ring of \textsf{finite representation type} if $A$ satisfies the following equivalent conditions.

\textbf{1)} $A$ is a right Artinian ring which has, up to isomorphism, only finite number of finitely generated indecomposable right modules.

\textbf{2)} $A$ is a left Artinian ring which has, up to isomorphism, only finite number of finitely generated indecomposable left modules.

\textbf{3)} $A$ is a right and left pure-semisimple ring.

The equivalence of the above conditions \textbf{1)--3)} is well-known; for example, see \cite{EisG71}, \cite{Aus74}, \cite{RinT75}, \cite{FulR75}. We note that the following question remains open: is it true that every right pure-semisimple ring is left pure-semisimple?

The first three main results of this paper are the following Theorems 1.3, 1.4 and 1.5. 

\textbf{Theorem 1.3.} Let $A$ be ring with minimum condition on principal right ideals. The following assertions are true.

\textbf{1.} The $\aleph_0$-distributive right $A$-modules coincide with the Artinian right $A$-modules. 

\textbf{2.} The $\aleph_0$-distributive left $A$-modules coincide with the Noetherian left $A$-modules.

\textbf{3.} The distributive left $A$-modules coincide with the completely cyclic left $A$-modules.

\textbf{4.} Distributive right $A$-modules coincide with completely cyclic Artinian right $A$-modules and, in particular, have composition series.

\textbf{5.} If the ring $A$ also satisfies the minimum condition on principal left ideals, then every $\aleph_0$-distributive right module is a Noetherian Artinian module.

\textbf{Theorem 1.4.} For a ring $A$, the following conditions are equivalent.

\textbf{1)} $A$ is a right $\aleph_0$-distributive ring with minimum condition on principal right ideals.

\textbf{2)} $A$ is a right $\aleph_0$-distributive ring with minimum condition on principal left ideals.

\textbf{3)} $A$ is a right Artinian ring.

\textbf{Theorem 1.5.} Every right module over the ring $A$ is a direct sum of $\aleph_0$-distributive modules if and only if
$A$ is a ring of finite representation type.\\
Consequently, every right module over the ring $A$ is a direct sum of distributive modules if and only if
$A$ is a ring of finite representation type and every indecomposable right $A$-module is distributive.

\textbf{Remark 1.6.} Two rings $A$ and $B$ are said to be \textsf{Morita equivalent} if the categories of right modules over these rings are equivalent. 

A ring $A$ is called a \textsf{right Koethe ring} if every right $A$-module is a direct sum of cyclic modules.

A ring $A$ is called a \textsf{right Kawada ring} if every ring,
which is Morita equivalent to the ring $A$, is a right Koethe ring.

We note that all Koethe rings and all Kawada rings are rings of finite representation type.

The fourth main result of this paper is Theorem 1.7 which uses the notion of the basis ring given in Remark 4.4(b) of Section 4 below. 

\textbf{Theorem 1.7.} For a ring $A$, the following conditions are equivalent.

\textbf{1)} Every right $A$-module is semidistributive.

\textbf{2)} $A$ is a Kawada ring with basis ring $B$ and every right $A$-module and every right $B$-module are direct sums of of completely cyclic distributive modules.

\textbf{3)} $A$ is a Kawada ring with basis ring $B$ and every indecomposable finitely generated right $B$-module is a completely cyclic module.

\textbf{Remark 1.8.} Let $A$ be a $F$-finite-dimensional algebra over a field $F$. Kawada's papers \cite{Kaw62}, \cite{Kaw65}, \cite{Kaw65} contain 19 rather complex conditions on the local idempotents of the basis ring $B$ of the algebra $A$, the fulfillment of which is equivalent to the fact that $A$ is a right Koethe ring. In \cite{Rin81}, this Kawada's result is analyzed and commented. To these 19 conditions, we can add the following condition 20: for every local idempotent $e$ of the algebra $B$, the module $eB_B$ is completely cyclic. Then one obtains a formal description of finite-dimensional algebras over a field which are rings, over which all modules are semidistributive. Of course, such description is not very useful.

\section{The proof of Theorems 1.3 and 1.4}\label{section2}

For the convenience of the reader, we present some well-known information.

We denote by $J(M)$, the Jacobson radical of the module $M$; we assume that $M=J(M)$ if the module $M$ does not have maximal submodules

A nonzero-left module is said to be \textsf{simple}\label{simmod} if it coincides with any its non-zero submodule. Direct sums of simple modules are called \textsf{semisimple}\label{semsimmod} modules. 
The semisimple modules coincide with non-zero modules, in which all submodules are direct summands. The right (left) semisimple rings coincide with rings which are are isomorphic to finite direct products of matrix rings over division rings.

A module $M$ is said to be \textsf{Artinian} (resp., \textsf{Noetherian}) if $M$ does not contain infinite properly descending (resp., properly ascending) chains of submodules. A module is a Noetherian if and only if all its submodules are finitely generated. Every right Artinian ring is right Noetherian. Every quasi-cyclic Abelian group is an Artinian $\mathbb{Z}$-module which is not Noetherian.
If $M$ is a semisimple module, then the module $M$ is finitely generated if and only if $M$ is a Noetherian Artinian module.

A module $M$ is called a \textsf{completely cyclic} module (resp., a \textsf{Bezout} module),\label{bezout} if all submodules (resp., all finitely generated submodules) of $M$ are cyclic.
Every uniserial module is a distributive Bezout module. The ring of integers is a commutative distributive non-uniserial Besout ring.
Any completely cyclic modules is Noetherian. The polynomial ring $\mathbb{Z}[x]$ over the ring of integers $\mathbb{Z}$ is a Noetherian $\mathbb{Z}[x]$-module which is not completely cyclic.

A non-zero module is said to be \textbf{semi-Artinian} if every its non-zero factor module contains a simple submodule. Every Artinian module is semi-Artinian. Direct sums of an infinite number of simple modules are semi-Artinian non-Artinian modules.

\textbf{Remark 2.1.} A module $M$ is said to be \textbf{semilocal} if $M$ is finitely generated and the factor module $M/J(M)$ with respect to its Jacobson radical is semisimple.

A ring $A$ is said to be \textsf{semilocal} if $A$ satisfies the following equivalent conditions.
\begin{itemize}
\item 
$A_A$ is a semilocal module.
\item 
$_AA$ is a semilocal module.
\item 
The factor ring $A/J(A)$ is isomorphic to a finite direct product of matrix rings over division rings.
\end{itemize} 
We note that every local ring is semilocal, and the direct product of two division rings is semilocal but is not local.

\label{!!!}

A ring $A$ is said to be \textsf{semiprimary} if $A$ is semilocal and its Jacobson radical $J(A)$ is nilpotent. Every right or left Artinian ring is semiprimary. If $\mathbb{R}$ and $\mathbb{Q}$ are the fields of real numbers and rational numbers, correspondingly, then the ring of matrices of the form
$$
\begin{pmatrix}
a&c\\
0&b
\end{pmatrix},\quad a,b\in\mathbb{Q},\; c\in\mathbb{Q}, 
$$
is semiprimary but is not a right or left Artinian ring. 

\textbf{Proposition 2.2.} Let $M$ be a $\aleph_0$-distributive module over the ring $A$ and the set of all non-isomorphic simple subfactors of the module $M$ finitely (for example, this is the case if the ring $A$ is semilocal). Then $M$ is a quotient finite-dimensional module.

\textbf{Proof.} Since all factor modules of $\aleph_0$-distributive modules are $\aleph_0$-distributive, it is sufficient to prove that $M$ is a finite-dimensional module. We assume that this is not true. Then $M$ contains a direct sum $\oplus_{k=1}^{\infty}X_k$ non-zero cyclic submodules $X_k$. Therefore, some factor module $\overline{M}$ of the module $M$ contains the direct sum $\oplus_{k=1}^{\infty}\overline{X_k}$ non-zero simple submodules $\overline{X_k}$. Since the set of all non-isomorphic simple subfactors of the module $M$ finitely, the set of simple modules $\overline{X_k}$ contains an infinite number of isomorphic modules. This This is a contradiction, since $M$ is an $\aleph_0$-distributive module.~$\square$

\label{!!!}

\textbf{Proposition 2.3.} Let $M$ be a $\aleph_0$-distributive semi-Artinian module over the ring $A$ and the set of all non-isomorphic simple subfactors of the module $M$ is finite (for example, this is the case if the ring $A$ is semilocal). Then $M$ is an Artinian module.

\textbf{Proof.} We assume the contrary. Then there exists an infinite properly descending chain $X_1\supset X_2\supset \ldots$ of submodules of the module $M$. We set $X=\cap_{k=1}^{\infty}X_k$, $\overline{M}=M/X$ and $\overline{X_k}=X_k/X$. A module $\overline{M}$ contains an infinite properly descending chain $\overline{X_1}\supset \overline{X_2}\supset \ldots$ of submodules which has the zero intersection. This is a contradiction, since the semi-Artinian module $\overline{M}$ is finite-dimensional, by Proposition 2.2; therefore, it is an essential extension of some finite direct sum $S$ of simple modules which is an Artinian module.~$\square$

A submodule $X$ of the module $M$ is said to be \textsf{small} in $M$ if $X+Y\ne M$ for any proper submodule $Y$ of the module $M$.

\textbf{Remark 2.4.} A ring $A$ satisfies minimum condition on principal right ideals if and only if the ring $A$ is semilocal and for any sequence $(x_k)_{k=1}^{\infty}\subset J(A)$, the element $x_1x_2\ldots x_n$ is equal to the zero for some $n$ depending on the sequence $(x_k)$) in general case; \cite{Bas60} see, also \cite[Theorem 11.6.3]{Kas82}. Such rings are also called \textsf{left perfect} rings. Under these conditions ,for every left $A$-module $M$, its submodule $J(A)M$ is small in $M$; see \cite[Theorem 11.5.7]{Kas82}.

\textbf{Remark 2.5.} A submodule $M$ of the module $X$ is said to be \textsf{fully invariant} submodule in $X$ if $\alpha (M)\subseteq M$ for every endomorphism $\alpha$ of the module $X$. A module $M$ is said to be \textsf{quasi-invariant} if all its maximal submodules are fully invariant in $M$.

The ideals of the ring $A$ coincide with fully invariant submodules in $A_A$ and with fully invariant submodules in ${}_AA$. Therefore, a ring is \textsf{right quasi-invariant} (resp., \textsf{left quasi-invariant}) if and only if all its maximal right (resp., left) ideals are ideals. It is clear that the ring $A$ is quasi-invariant right if and only if the factor ring $A/J(A)$ is quasi-invariant right. Consequently, the factor ring $A/J(A)$ is a direct product of division rings, then the ring $A$ is (right and left) quasi-invariant.

\textbf{Proposition 2.6; \cite{Ste74}.} Let $M$ be a distributive module. The following assertions are true.

\textbf{1.} $M$ does not have subfactors which are direct sums of two isomorphic non-zero modules.

\textbf{2.} A module $M$ is quasi-invariant.

\textbf{3.} Every submodule of the module $M$ which is a finite direct sum of cyclic of modules, is cyclic. 

\textbf{Proposition 2.7.} If $A$ is a semilocal ring, then every finitely generated distributive right $A$-module $M$ is a cyclic module, i.e., all distributive right $A$-modules are Bezout modules.

\textbf{Proof.} Since $M/MJ(A)$ is a distributive semisimple module over the semisimple Artinian ring $A/J(A)$, then the factor module $M/MJ(A)$ is an Artinian semisimple module, by Proposition 2.3. Therefore, distributive module $M/MJ(A)$ is a finite direct sum of simple of modules. By Proposition 2.6(3), the module $M/MJ(A)$ is cyclic and the submodule $MJ(A)$ is small in $M$, by the Nakayama lemma. Therefore, the module $M$ is cyclic.~$\square$

\textbf{Proposition 2.8.} If $A$ is a ring with minimum condition on principal right ideals, then every $\aleph_0$-distributive right $A$-module $M$ is a Artinian module.

\textbf{Proof.} Since $A$ is a ring with minimum condition on principal right ideals, it is directly verified that $M$ is a semi-Artinian module. It follows from Proposition 2.3 that the module $M$ is Artinian.~$\square$

\textbf{Proposition 2.9.} Let $A$ be ring with minimum condition on principal left ideals (i.e., the ring $A$ is right perfect).

\textbf{1.} Every $\aleph_0$-distributive right $A$-module $M$ is a Noetherian module. 

\textbf{2.} Every distributive right $A$-module $M$ is a completely cyclic Noetherian module. 

\textbf{Proof.} \textbf{1.} Since all submodules $\aleph_0$-distributive of modules are $\aleph_0$-distributive, it is sufficient to prove that module $M$ is finitely generated. Since $M/MJ(A)$ is an $\aleph_0$-distributive semisimple module over the semisimple Artinian ring $A/J(A)$, it follows from Proposition 2.3 that $M/MJ(A)$ is an Artinian semisimple module; thefore, it is finitely generated. In addition, it folows from Remark 2.4 that the submodule $MJ(A)$ is small in $M$, since $M$ is a right module over the right perfect ring $A$. Therefore, the module $M$ is finitely generated.

\textbf{2.} Let $X$ be a submodule of the module $M$. It follows from \textbf{1} that the module $X$ is finitely generated. Since the ring $A$ is semilocal, it follows from Proposition 2.7 that the module $X$ is cyclic.~$\square$

\textbf{Remark 2.10.} It follows from Propositions 2.8 and 2.9(1) that every $\aleph_0$-distributive right module over the ring with minimum condition on principal left ideals and with minimum condition on principal left ideals for principal right ideals is an  Artinian Noetherian module. In particular, every $\aleph_0$-distributive right module over a semiprimary ring is an Artinian Noetherian module.

\textbf{Remark 2.11.} \textbf{The completion of the proof of Theorems 1.3 and 1.4.}

\textbf{Proof.} Item \textbf{1} of Theorem 1.3 follows from Remark 1.1 and Proposition 2.8.

Item \textbf{2} of Theorem 1.3 follows from Remark 1.1 and Proposition 2.9(1).

Item \textbf{3} of Theorem 1.3 follows from of Theorem 1.3(2) and Proposition 2.7.

Item \textbf{4} of Theorem 1.3 follows from Theorem 1.3(1) and Proposition 2.7.

Item \textbf{5} of Theorem 1.3 follows from Theorem 1.3(1) and Theorem 1.3(2).

In Theorem 1.4, the implication 3)\,$\Rightarrow$\,1) follows from Remark 1.1, the implication 3)\,$\Rightarrow$\,2) follows from Remark 1.1 and the property that every right Artinian ring is right and left perfect, and the implication 1)\,$\Rightarrow$\,3) follows from Proposition 2.8. 

We prove the implication 2)\,$\Rightarrow$\,3) of Theorem 1.4. By Proposition 2.9(1), the ring $A$ is right Noetherian. In addition, the right perfect ring $A$ is semilocal and its Jacobson radical is a nil-ideal. Since any nil-ideal of the right Noetherian ring $A$ is nilpotent, the ring $A$ is semiprimary. It is clear that the right Noetherian semiprimary ring $A$ is right Artinian.~$\square$

\section{The proof of Theorem 1.5}\label{section3}

\textbf{Remark 3.1.} A module $M$ is said to be \textsf{injective with respect to the module}\label{relinj} $X$ or \textsf{$X$-injective} if for any submodule $X_1$ in $X$ every homomorphism $X_1\to M$ is extended to homomorphism $X\to M$. A module $M$ over the ring $A$ is said to be \textsf{injective}\label{injmod} if $M$ is injective with respect to any $A$-module. For example, over finite direct products of matrix rings over division rings, all modules are injective. In addition, an Abelian group $M$ is an injective module over the ring of integers $\mathbb{Z}$ if and only if $M$ is a divisible Abelian group, i.e., $M$ is a direct sum of groups, isomorphic to the additive group $\mathbb{Q}$ of rational numbers and quasi-cyclic groups $\mathbb{Z}(p^{\infty})$.

A module $Q$ is called an \textsf{essential extension} of its submodule $M$ if $M\cap X\ne 0$ for any non-zero (cyclic) submodule $X$ of the module $Q$. In this case, $M$ is called an \textsf{essential submodule} in $Q$. Every module $M$ is an essential submodule of some injective module $Q$ which is called an \textsf{injective hull}\label{injhu} of the module $M$ and the injective hull $Q$ is unique up to isomorphism. For example, if $\mathbb{Z}$ is the ring of integers, then the additive group $\mathbb{Q}$ is an injective hull of the module $\mathbb{Z}_{\mathbb{Z}}$.

\textbf{Proposition 3.2.} For any ring $A$, there exists a cardinal number $\aleph$ such that the cardinality of any non-zero $\aleph_0$-distributive right $A$-module $M$ does not exceed $\aleph$.

\textbf{Proof.} For any cyclic right $A$-module $X$, there exists a right ideal $B$ of the ring $A$ with $X\cong A_A/B$. Therefore, we can consider the \textbf{set} $\{X_i\}_{i\in I}$ of all pairwise non-isomorphic cyclic right $A$-modules. Let $Y_i$ be the direct sum of the countable set of isomorphic copies of the module $X_i$, $i\in I$. Let $Q$ be the injective hull of the module $\oplus_{i\in I}Y_i$ and let $\aleph$ be the the cardinality of the module $Q$. With the use of the Zorn lemma, we can verify that the module $M$ is an essential extension of some submodule $\oplus_{j\in J}M_j$, where all $M_j$ are non-zero cyclic modules. Since the module $M$ is $\aleph_0$-distributive, there exists a monomorphism $f\colon \oplus_{j\in J}M_j\to \oplus_{i\in I}Y_i$. Since the module $Q$ is injective, the monomorphism $f$ is extended to a homomorphism $g\colon M\to Q$. In addition, $(\oplus_{j\in J}M_j)\cap \text{Ker } g(\oplus_{j\in J}M_j)\cap \text{Ker } f=0$. Therefore, $g$ is a monomorphism, since $M$ is an essential extension of the module $\oplus_{j\in J}M_j$. Then module $M$ is isomorphic to a submodule of the module $Q$. Therefore, the cardinality of the module $M$ does not exceed $\aleph$.~$\square$ 

\textbf{Proposition 3.3; \cite{Cha60}, also see \cite[Theorem 20.23]{Fai76}.}\\ If there exists a cardinal number $\aleph$ such that any right $A$-module is a direct sum of $\aleph_0$-generated modules, then the ring $A$ is right Artinian.

\textbf{Remark 3.4.} \textbf{The completion of the proof of Theorem 1.5.}

\textbf{Proof.} If $A$ is a ring of finite representation type, then every $A$-module is a direct sum of indecomposable finitely generated modules; these modules are Artinian and, in particular, $\aleph_0$-distributive.

Now let every right module over the ring $A$ be a direct sum of $\aleph_0$-distributive modules. It follows from Theorem 1.3(5) that every $\aleph_0$-distributive right $A$-module is a Noetherian Artinian module. Therefore, $A$ is a ring of finite representation type.~$\square$ 

\textbf{Proposition 3.5; \cite[Theorem 3.26]{Tug98}.} If $A$ is a right quasi-invariant ring, then every right Bezout $A$-module is distributive.

\section{The proof of Theorem 1.7}\label{section4}

In Remarks 4.1, 4.2 and 4.3 below, we give some well-known information on semiperfect rings.

\textbf{Remark 4.1.} A module $M$ is said to be \textbf{local},\label{locMod} if $M$ is cyclic and the factor module $M/J(M)$ with respect to its Jacobson radical is simple.

A ring $A$ is said to be \textsf{local},\label{locRin} if $A$ satisfies the following equivalent conditions.

\textbf{1)} $A_A$ is a local module.

\textbf{2)} $_AA$ is a local module.

\textbf{3)} The factor ring $A/J(A)$ is a division ring.

\textbf{4)} For any element $a\in A$, at least one of the elements $a,1-a$ is invertible.

A non-zero idempotent $e$ of the ring $A$ is said to be \textsf{local}\label{locRin} if the ring $eAe$ is local.

\textbf{Remark 4.2.} Let $M$ be a distributive module over the ring $A$. If all simple subfactors of the module $M$ are isomorphic to each other (for example, it is the case if the ring $A$ is local), then $M$ is a uniserial module. 

\textbf{Proof.} Indeed, if $M$ is not uniserial, then $M$ contains two cyclic submodules $X$ and $Y$ which are incomparable with respect to inclusion. Then the submodule $X+Y$ has a factor module $S\oplus T$, where the modules $S$ and $T$ are simple. By assumption, the simple modules $S$ and $T$ are isomorphic to each other; this contradicts to Proposition 2.6(1).~$\square$

\textbf{Remark 4.3.} A ring $A$ is said to be \textsf{semiperfect} if $A$ satisfies the following equivalent conditions.

\textbf{1)} The ring $A$ is semilocal and all idempotents of the factor ring $A/J(A)$ are lifted to idempotents of the ring $A$.

\textbf{2)} $A_A$ is a (finite) direct sum of of local modules.

\textbf{3)} $_AA$ is a (finite) direct sum of of local modules.
 
\textbf{4)} There exists a decomposition $1=e_1+\ldots+e_n$ of the identity element of the ring $A$ into a sum of local pairwise orthogonal idempotents $e_1,\ldots,e_n$.

We note that if $p$ and $q$ are two distinct prime integers and the ring $A$ consists of all rational numbers such that which denominators are not divided by $p$ or $q$ in irreducible form, then $A$ is a semilocal non-semiperfect ring.

All semiprimary rings and serial right or left rings are semiperfect. For any division ring $F$, the ring $F[[x,y]]$ of formal power series in two independent variables $x$ and $y$ is semiperfect (in fact, it is local) but it is neither a semiprimary ring, nor a right serial ring, nor a left serial ring.

\textbf{Remark 4.4.} Let $A$ be a semiperfect ring. 

\textbf{a.} The decomposition $1=e_1+\ldots+e_n$ from Remark 4.3(4) is unique up to isomorphism, i.e., for any other decomposition $1=f_1+\ldots+f_m$ of the identity element of the semiperfect ring $A$ into the sum of local pairwise orthogonal idempotents $f_1,\ldots,f_m$, we have that $m=n$ and there exists a permutation $\tau$ of the set $\{e_1,\ldots,e_n\}$ such that $e_iAe_i\cong f_{\tau(i)}Af_{\tau(i)}$, $e_iA\cong f_{\tau(i)}A$ and $Ae_i\cong Af_{\tau(i)}$ for all $i\in \{1,\ldots,n\}$.

In addition, the sets $\{e_1,\ldots,e_n\}$ and $\{f_1,\ldots,f_n\}$ contain some subsets $\{e_1',\ldots,e_k'\}$ and $\{f_1',\ldots,f_{\ell}'\}$, correspondingly, such that $e_i'A\not\cong e_j'A$ and $f_i'A\not\cong f_j'A$ for $i\ne j$ (equivalently, $Ae_i'\not\cong Ae_j'$ and $Af_i'\not\cong Af_j'$ for $i\ne j$), and $k=\ell$ and the categories of all right $(e_1'+\ldots+e_k')A(e_1'+\ldots+e_k')$-modules and of all right $(f_1'+\ldots+f_k')A(f_1'+\ldots+f_k')$-modules are equivalent to the category of all right $A$-modules (a similar assertion is true for the corresponding categories of left modules).

\textbf{b.} The idempotents $e'=e_1'+\ldots+e_k'$ and $f'=f_1'+\ldots+f_k'$, which are mentioned in \textbf{a}, are called \textsf{basis idempotents} of the ring $A$, and the rings $e'Ae'$ and $f'Af'$ are isomorphic to each other and are called \textsf{basis rings} of the ring $A$. If the basis idempotent $e'$ is equal to of the identity element of the ring $A$ (i.e., $A$ coincides with its basis ring), then $A$ is called a \textbf{basis} (or \textsf{self-basis}) ring. 

\textbf{Remark 4.5.} Let $B$ be the basis ring of the semiperfect ring $A$. The following assertions \textbf{a}, \textbf{b} and \textbf{c} are well-known. The assertion of \textbf{d} follows from the assertions of \textbf{c}, Proposition 3.5 and Proposition 2.9(2). The assertion of \textbf{e} follows from \textbf{d}.

\textbf{a.} The factor ring $B/J(B)$ is isomorphic to a direct product of division rings. 

\textbf{b.} $B$ is a basis, semiperfect, quasi-invariant ring. 

\textbf{c.} The rings $A$ and $B$ are Morita equivalent, i.e., there exists an equivalence functor $\mathcal{T}$ between the category of all right $A$-modules $\text{Mod }A$ and the category of all right $B$-modules $\text{Mod }B$ and the ring $A$ is right (resp., left) perfect if and only if the ring $B$ is right (resp., left) perfect. 
In particular, distributivity (resp., semidistributivity) of some right $A$-module $M$ if and only if the right $B$-module $\mathcal{T}(M)$ is distributive (resp., semidistributive).

\textbf{d.} Distributivity (resp., semidistributivity) of some right $A$-module $M$ is equivalent to the property that the right $B$-module $\mathcal{T}(M)$ is a Bezout module (resp., direct sum of Bezout modules). In addition, if the ring $A$ (and, consequently, the ring $B$) is right perfect, then distributivity (resp., semidistributivity) of the right $A$-module $M$ is equivalent to the property that the right $B$-module $\mathcal{T}(M)$ is completely cyclic (resp., is a direct sum of of completely cyclic modules).

\textbf{e.} Every right module over the ring $A$ is a direct sum of distributive modules if and only if
$A$ is a Artinian ring of finite representation type with basis ring $B$ and every indecomposable right $B$-module is completely cyclic.

\textbf{Remark 4.6.} \textbf{The completion of the proof of Theorem 1.7.}

\textbf{Proof.} The implication 2)\,$\Rightarrow$\,3) is trivial.

The implication 3)\,$\Rightarrow$\,1) follows from Remark 4.5(e).

1)\,$\Rightarrow$\,2). Let every right $A$-module be semidistributive. Since the property of semidistributivity of all right modules is preserved under Morita equivalence, it follows from Remark 4.5(e) that $A$ is an Artinian ring of finite representation type with basis ring $B$ and every indecomposable right $B$-module is completely cyclic. In particular, $A$ is a right Koethe ring. Since semidistributivity property of all right modules is preserved under Morita equivalence, every ring, which is Morita equivalent to the ring $A$, is a Koethe ring, i.e., the rings $A$ and $B$ are Kawada rings.~$\square$


\begin{thebibliography}{99}

\bibitem{Aus74} Auslander M. 1974. Representation theory of Artin algebras II. Comm. Algebra 1: 269-310.

\bibitem{Cha60} Chase S.U. Direct products of modules. -- Trans. Amer. Math. Soc. -- 1960. -- Vol.~97. -- P.457-473.

\bibitem{EisG71} Eisenbud D., Griffith Ph. Serial rings. -- J.~Algebra. -- 1971. -- Vol.~17. -- P.~389-400.

\bibitem{Fai76} Faith C. Algebra II. -- Springer, Berlin -- New York. -- 1976.

\bibitem{FulR75} Fuller K.R., Reiten I. Note on rings of finite representation type and decompositions of modules. -- Proc. Amer. Math. Soc. -- 1975. -- Vol.~50. -- P.~92-94.

\bibitem{Kas82} Kasch F. Modules and Rings. -- Academic Press, London? 1982.

\bibitem{Kaw62} Kawada Y. On Koethe’s problem concerning algebras for which every indecomposable module is cyclicI, Sci. Rep. Tokyo Kyoiku Daigaku Sect.A 7 (1962) 154–230.

\bibitem{Kaw63} Kawada Y., On Koethe’s problem concerning algebras for which every indecomposable module is cyclic II, Sci. Rep. Tokyo Kyoiku Daigaku Sect.A 8 (1963) 1–62.

\bibitem{Kaw65} Kawada Y. On Koethe’s problem concerning algebras for which every indecomposable module is cyclic III, Sci. Rep. Tokyo Kyoiku Daigaku Sect.A 8 (1965) 165–250.

\bibitem{Rin81} Ringel C.M. Kawada’s theorem. Lecture Notes Math. -- 1981. -- Vol.~874. -- Springer, Berlin -- New York. -- P.431-447.

\bibitem{RinT75} Ringel C.M., Tachikawa H. QF-3 rings. -- J. Reine Angew. Math. -- 1975. -- Vol.~272. -- P.~49-72.

\bibitem{Ste74} Stephenson W. Modules whose lattice of submodules is distributive. -- Proc. London Math. Soc. -- 1974. -- Vol.~28, no.~2. -- P.~291-310.

\bibitem{Tug98} Tuganbaev A.A. Semidistributive Modules and Rings. -- Dordrecht -- Boston -- London, Kluwer Academic Publishers, 1998.

\end{thebibliography}
\end{document}